\documentclass[11pt]{amsart}
\usepackage{amssymb,latexsym,comment,url,amsmath}

\usepackage[T1]{fontenc}
\usepackage{rotating,graphicx}

\usepackage{diagbox}

\usepackage{longtable}

\newcommand{\EE}{{\mathcal E}}

\newcommand{\Fp}{{\mathbb{F}_p}}

\newcommand{\PP}{{\mathbb P}}
\newcommand{\Q}{{\mathbb Q}}

\newcommand{\ord}{\operatorname{ord}}
\newcommand{\rank}{\operatorname{rank}}

\newcommand{\Z}{{\mathbb Z}}

\newenvironment{Proof}{\par\noindent{\sc Proof:}}%
                      {\hspace*{\fill}\nobreak$\Box$\par\medskip}
                       {\hspace*{\fill}\nobreak$\Box$\par\medskip}

\newtheorem{Proposition}{Proposition}[section]
\newtheorem{Theorem}[Proposition]{Theorem}
\newtheorem{Lemma}[Proposition]{Lemma}

\theoremstyle{definition}
\newtheorem{Definition}[Proposition]{Definition}
\newtheorem{Remark}[Proposition]{Remark}
 
\renewcommand{\baselinestretch}{1.1}

\begin{document}

\title[Rational points on quadratic elliptic surfaces]%
{Rational points on quadratic elliptic surfaces}

\author[M. Sadek]%
{Mohammad~Sadek}
\address{Faculty of Engineering and Natural Sciences, Sabanc{\i} University, Tuzla, \.{I}stanbul, 34956 Turkey}
\email{mohammad.sadek@sabanciuniv.edu}

\begin{abstract}
{We consider elliptic surfaces whose coefficients are degree $2$ polynomials in a variable $T$. It was recently shown, \cite{Kollar}, that for infinitely many rational values of $T$ the resulting elliptic curves have rank at least $1$. In this article, we prove that the Mordell-Weil rank of each such elliptic surface is at most $6$ over $\Q$. In fact, we show that the Mordell-Weil rank of these elliptic surfaces is controlled by the number of zeros of a certain polynomial over $\Q$.}
\end{abstract}

\maketitle

\let\thefootnote\relax\footnotetext{ 
The datasets generated during and/or analysed during the current study are available from the corresponding author on reasonable request.}

\section{Introduction}
Let $\EE$ be a nonconstant elliptic surface defined by the following equation over $\Q$
\begin{eqnarray} \label{eq*}y^2 = a_3(T)x^3 + a_2(T)x^2 + a_1(T)x + a_0(T),\end{eqnarray}
where $a_i(T)\in\Q[T]$. The elliptic surface $\EE$ can be thought of as a family of elliptic curves parametrized by $T\in\Q$. The Lang-N\'{e}ron theorem asserts that the abelian group $\EE(\Q(T))$ of rational points is finitely generated. In addition, for all but finitely many rational values of $T=t$, specialization yields an elliptic curve $\EE_t$. The Mordell-Weil theorem states that the abelian group $\EE_t(\Q)$ is also finitely generated. Thanks to Silverman's specialization theorem, one knows that for all but finitely many values $t\in\Q$, the Mordell-Weil rank of $\EE_t(\Q)$ is at least the Mordell-Weil rank of $\EE(\Q(T))$.

The construction of elliptic surfaces over $\Q$ with large Mordell-Weil rank has been an interesting target for many researchers. One of the main reasons is that among specializations of these elliptic surfaces one may find elliptic curves with large Mordell-Weil rank over $\Q$. In \cite{Mest1, Mest2}, Mestre gave examples of elliptic surfaces of Mordell-Weil rank $\ge11,12$ over $\Q(T)$. Nagao used ideas of Mestre to search for elliptic surfaces of larger Mordell-Weil rank over $\Q(T)$, see \cite{Nagao2}. The arithmetic of these elliptic surfaces was employed to find elliptic curves of Mordell-Weil ranks $\ge 19,20, 21, 22$ in \cite{Fer1}, \cite{Nagao1}, \cite{Nagao3}, and \cite{Fer2}, respectively.

The following question is frequently posed: can we find rational values for $t$ such that the specialized elliptic curves $\EE_t$ satisfy that
$\rank(\EE_t(\Q)) \ge \rank(\EE(\Q(T))) + 1$? If the answer is positive, then are there infinitely many such rational values for $t$?
 If the answer again is in the affirmative, then is the density (appropriately defined) of such rational values for $t$ positive?

Elliptic surfaces defined by equation (\ref{eq*}) over $\Q$ where the polynomials $a_i(t)$ are of degree at most $1$ were studied in \cite{Mackall}. It was proven that the Mordell-Weil rank of such an elliptic surface is at most $3$,  see \cite[Lemma 3.1]{Mackall}. Explicit families of such elliptic surfaces with Mordell-Weil rank $0,1,2,$ or $3$ can be found in \cite{Fer,Mackall}.

In this article, we consider the elliptic surface $\EE$ defined by
$y^2 = a_3(T)x^3 + a_2(T)x^2 + a_1(T)x + a_0(T)$,
 where $a_i (T) \in \Q[T]$ are polynomials of degree at most $2$. We assume moreover it is a nontrivial elliptic surface, i.e., at least two of the curves $\EE_t$ are smooth, elliptic and not isomorphic to each other over $\Q$. This quadratic family of elliptic curves was studied in \cite{Kollar}. The authors prove that when the whole family is viewed as a single algebraic surface in $\mathbb{A}^3_{xyt}$, then the resulting surface is unirational over $\Q$. As a consequence, it is proved that there are infinitely many rational values $T = t$ such that $\rank(\EE_t(\Q)) > 0$. In fact, it is conjectured that such rational values constitute a positive
proportion of all rational numbers. A key ingredient to reach the aforementioned results is the observation
$$a_3(T)x^3 + a_2(T)x^2 + a_1(T)x + a_0(T)=A(x)T^2+B(x)T+C(x)$$
where $A(x),B(x),C(x)\in\Q[x]$ are of degree at most $3$. In particular, $\EE$ is birational to a conic bundle with at most $7$ singular fibers.

Nagao conjectured that the Mordell-Weil rank of an elliptic surface can be recovered from local information of the elliptic surface modulo primes, see \cite{Nagao}. This conjecture was proved unconditionally for rational elliptic surfaces, and conditionally on Sato-Tate conjecture for other elliptic surfaces by Rosen and Silverman, see \cite{Rosen}. We show that the quadratic elliptic surface
$$\EE : y^2 = a_3(T)x^3 + a_2(T)x^2 + a_1(T)x + a_0(T)=A(x)T^2+B(x)T+C(x),$$
are rational elliptic surfaces under some mild conditions. 

Making use of the fact that the aformentioned quadratic elliptic surfaces are rational enables us to evaluate their Mordell-Weil rank. In fact, we show that the Mordell-Weil rank depends on the factorization behavior of the polynomial $D(x):=B(x)^2-4A(x)C(x)$ in $\Q[x]$ where the latter polynomial $D(x)$ is of degree at most $6$. In particular, we show that the largest possible value for the Mordell-Weil rank $r_{\EE}$ of $\EE$ is  $6$. Moreover, this maximal rank occurs if and only if $D(x)$ splits completely into distinct linear factors over $\Q$ and the value of $A(x)$ at each of the simple zeros of $D(x)$ is a nonzero rational square. In addition, we show the existence of quadratic elliptic surfaces with Mordell-Weil rank $r$ where $0 \le r\le   6$. The aforementioned result asserts that for all but finitely many rational values of $t$, the Mordell-Weil rank of the specialized elliptic curve $\EE_t$ is at least $r_{\EE}$. This provides an improvement on the results of \cite{Kollar} mentioned above.

It is worth mentioning that elliptic surfaces defined by equation (\ref{eq*}) over $\Q$, where at least one of the polynomials $a_i(t)$ is of degree at least $3$, are not in general rational elliptic surfaces. However, assuming Sato-Tate conjecture for these surfaces, Nagao's conjecture on their Mordell-Weil ranks holds true, see \cite{Rosen}. Yet, since these elliptic surfaces are not birational to conic bundles, the techniques used in this work can not be applied to evaluate their Mordell-Weil ranks. 

After the submission of this article, it came to the attention of the author that the authors of \cite{Bat} have been using similar techniques to find expressions for the rank of some quadratic elliptic surfaces over $\Q(T )$ depending on the factorization and other properties of certain polynomials. In fact, the authors treat elliptic surfaces over $\Q$ whose Weierstrass equations can be written as $$y^2=A(x)T^2+B(x)T+C(x), $$  where $A(x),B(x)\in\Q[x]$ are of degree at most $2$ and $ C(x)\in\Q[x]$ is monic of degree $3$, and one of $A(X)$ or $B(X)$ is not the zero polynomial. The authors then introduce conditions under which these elliptic surfaces possess positive Mordell-Weil ranks. This is approached by considering several subcases based on whether $A(x)$ is a nonconstant polynomial, zero, a nonzero rational square, or a rational number which is not a rational square.  

\subsection*{Acknowledgments} The author would like to thank the anonymous referees for several comments and suggestions. This work is partially supported by BAGEP Award of the Science Academy, Turkey.

\section{Elliptic surfaces}
 \label{sec:elliptic surfaces}
Let $\EE\to \PP^1$ be a (non-split) elliptic surface over $\Q$ described by a Weierstrass equation
$$\EE : y^2 = x^3 +A(T)x+B(T), \,A,B \in \Z[T],\, \Delta(T) = 4A(T)^3 +27B(T)^2\ne 0.$$
By evaluating the polynomials $A(T)$ and $B(T)$ at integer values for $T$, we obtain elliptic curves over $\Q$.

 Silverman’s Specialization Theorem implies that, for all but finitely many values $t\in\Z$, the Mordell-Weil rank of the fiber $\EE_t$ over $\Q$ is at least that of the curve $\EE$ over $\Q(T)$.

Let $t\in\Z$ and let $p$ be a prime. If $p|\Delta(t)$, then we let $a_p(\EE_t)=0$. Otherwise, we let $a_p(\EE_t)$ be the $p$-th coefficient of the $L$-series of the fiber $\EE_t$.  We set
$$A_p(\EE) =\frac{1}{p}\sum_{t=1}^pa_p(\EE_t)$$ to be the average of the $a_p$'s over the fibers.

Nagao \cite{Nagao} conjectured that
$$\lim_{X\to\infty} S(X,\EE) = \rank\EE(\Q(T)),\textrm{  where }S(X,\EE) =\frac{1}{X} \sum_{p\le X }  -A_p(\EE)\log p.$$
In fact, he proved the conjecture for several elliptic surfaces. In \cite{Rosen}, it was proved, among other results, that Nagao's conjecture is true for rational elliptic surfaces, namely, if $\EE$ is birational to $\PP^2$. In terms of a minimal Weierstrass equation $y^2 = x^3 + A(T )x + B(T )$ for $\EE$ over $\Q(T )$, an elliptic surface is rational if and only if one of the following two conditions holds:
\begin{itemize}
\item[(i)] $0 < \max\{3\deg A,2\deg B\} < 12$.
\item[(ii)] $3\deg A = 2\deg B = 12$ and $\ord_{T=0} T^{12}\Delta(T^{-1}) = 0$,
\end{itemize}
see \cite[Remark 1.3.1]{Rosen}

\section{Quadratic families of elliptic curves}
\label{sec:quadratic families}
Let $\EE$ be an elliptic surface defined over $\Q$ by the following equation
$$y^2 = a_3(T)x^3 + a_2(T)x^2 + a_1(T)x + a_0(T)$$
where $a_i \in\Z[T]$, $\max_i \deg a_i(T) = 2$, $i = 0,1,2,3$.

We only consider non-split elliptic surfaces. Therefore, at least two of the curves $\EE_t$ are smooth, elliptic and not isomorphic to each other over $\Q$; $a_3(T)$ is not identically $0$; and not all the $a_i(T)$ are constant multiples of the same square $(a T + b)^2$.

One uses the transformation $x\mapsto x/a_3(T)$ and $y\mapsto  y/a_3(T)$ to see that $\EE$ is described by $y^2 = x^3 +a_2(T )\,x^2 +a_1(T )a_3(T )\,x+a_0(T )a_3(T )^2$. Then one sees that the transformation $x \mapsto x-a_2(T )/3$ allows $\EE$ to be described by
$$y^2=x^3+(a_1(T)a_3(T)-a_2(T)^2/3)x+a_0(T)a_3(T)^2-a_1(T)a_2(T)a_3(T)/3+2a_2(T)^3/27.$$
The discriminant $\Delta(T)$ of $\EE$ is
$$a_3(T)^2\left[-a_1(T)^2a_2(T)^2+4a_1(T)^3a_3(T)-18a_0(T)a_1(T)a_2(T)a_3(T)+a_0(T)\left(4a_2(T)^3+27a_0(T)a_3(T)^2\right)\right].$$

\begin{Lemma}
Let $\EE$ be an elliptic surface defined over $\Q$ by the following equation $$y^2 = a_3(T)x^3 + a_2(T)x^2 + a_1(T)x + a_0(T)$$
where $a_i(T)\in \Z[T]$ and $\deg a_i(T) \le 2$.

If the fiber $\EE_{\infty}$ is an elliptic curve, then $\EE$ is a rational elliptic surface. In particular, $\lim_{X\to\infty}S(X,\EE) =\rank(\EE (\Q(T )))$.
\end{Lemma}
\begin{Proof}
We will show that the condition $\EE_{\infty}$ is an elliptic curve implies that $\EE$ is an elliptic surface that satisfies $\ord_{T =0} T^{12} \Delta(T^{-1} ) = 0$, hence $\EE$ is a rational elliptic surface, see \S \ref{sec:elliptic surfaces}. The formula for $\Delta(T )$ implies that $T^{12}\Delta(T^{-1})$ at $T = 0$ is given by
$$\Delta=a_3^2\left(-a_1^2a_2^2+4a_1^3a_3-18a_0a_1a_2a_3+a_0\left(4a_2^3+27a_0a_3^2\right)\right)$$
where $a_i$ is the coefficient of $T^2$ in $a_i(T)$. But $\Delta$ is the discriminant of the cubic curve
$$\EE_{\infty}:y^2=x^3+(a_1a_3-a_2^2/3)x+a_0a_3^2-a_1a_2a_3/3+2a_2^3/27.$$
If $\EE_{\infty}$ is an elliptic curve, then $\Delta\ne 0$, hence the result.
\end{Proof}
If $\EE_{\infty}$ is an elliptic curve, then $\Delta\ne0$, and so $a_3\ne 0$. It is easy to see that 
$$a_3(T)x^3 + a_2(T)x^2 + a_1(T)x + a_0(T)=A(x)T^2+B(x)T+C(x)$$
where $\deg A = 3$ and $\deg B, \deg C \le 3$.  This holds as the leading coefficient of $A(x)$ is $a_3$. 

The elliptic surface $\EE$ can be viewed as the family of conics
$$y^2 =A(x)T^2 +B(x)T +C(x)$$
over $\Q(x)$. The singular conics are exactly those for which $x$ is a root of the discriminant $B(x)^2 - 4A(x)C(x)$. There are in general $6$, not necessarily distinct, such values for $x$. One may see that, after a birational transformation, the conic at infinity is isomorphic to the conic described by $a_3' (S, T ) = 0$ in $\PP^2_{STW}$,  where $a_3' (S, T ) :=  S^2a_3(T/S)$ is the homogenization of $a_3(T)$. If $a_3(T)$ is irreducible in $\Q[T]$, then the conic at infinity is singular, hence there are $\le 7$ singular conics in the family. If $a_3(T)$ is reducible, then the conic at infinity is reducible and we can contract either of its irreducible components. In the latter case, the conic bundle has at most $6$ singular conics, see \cite[\S 0]{Kollar}. Therefore, from now on we assume that $a_3(T )$ is irreducible of degree $2$.

\begin{Definition}
 Let $\EE$ be a non-split elliptic surface defined over $\Q$ by the equation $$y^2 = a_3(T)x^3 + a_2(T)x^2 + a_1(T)x + a_0(T)$$
where $a_i(T) \in \Z[T]$, $\deg a_3(T) = 2$, $\deg a_i(T) \le 2$, $ i = 0,1,2$, and $a_3(T)$ is irreducible. Assume moreover that $\EE_{\infty}$ is an elliptic curve. Then $\EE$ will be called a {\em quadratic elliptic surface}.
\end{Definition}

\section{Evaluation of some character sums}
\label{sec:character sums}
Let $p$ be a prime. Throughout the paper, $\Fp$ will be the finite field with $p$ elements, and $\phi_p$ will be the unique quadratic character modulo $p$.

One will need the following standard lemma to evaluate character sums of quadratic polynomials, see for example \cite[Lemma A.2]{Arms}.

\begin{Lemma}
\label{Lem:Legendre}
Let $p>2$ be a prime. Let $a,b,c\in\Fp$ be such that either $a$ or $b$ are nonzero. One then has
\begin{eqnarray*}
\sum_{t=0}^{p-1}\phi_p(at^2+bt+c)= \begin{cases} 
      (p-1)\phi_p(a) & \textrm{ if }p|(b^2-4ac) \\
      -\phi_p(a) &\textrm{ otherwise}  
   \end{cases}
\end{eqnarray*}
\end{Lemma}

\begin{Lemma}
\label{Lem:symbols of cubics}
Let $f(x)\in\Fp[x]$ be a cubic polynomial. One has
\begin{eqnarray*}
\sum_{x\in\Fp}\phi_p(f(x))= \begin{cases} 
     -a_p(E) & \textrm{ if $f(x)$ has no repeated root} \\
       -\phi_p(-f_1^{-1}f_2f_3+f_4)&\textrm{ if }f(x)=(f_1x+f_2)^2(f_3x+f_4),\;f_1f_3\in\mathbb{F}_p^{\times}  \\
       0& \textrm{ if }f(x)=(f_1x+f_2)^3,\;f_1\in\mathbb{F}_p^{\times}
   \end{cases}
\end{eqnarray*}
where $a_p(E)$ is the trace of Frobenius of the elliptic curve $E$ defined by $y^2=f(x)$ over $\Fp$.
\end{Lemma}
\begin{Proof}
If $f(x)$ has no repeated root, then one has $a_p(E) =-\sum_{x\in\Fp}
\phi_p(f(x))$. 

One now assumes that $f(x)$ has a repeated root. If $f(x) = (f_1 x+f_2 )^2(f_3x+f_4 )$, with $f_1f_3\in\mathbb{F}_p^{\times}$, then 
\begin{eqnarray*}
\sum_{x\in\Fp}\phi_p(f(x))&=&\sum_{x\in\Fp}\phi_p^2(f_1x+f_2)\phi_p(f_3x+f_4)=\sum_{x\ne -f_1^{-1}f_2}\phi_p(f_3x+f_4)\\
&=&\sum_{x\in\Fp}\phi_p(f_3x+f_4)-\phi_p(-f_1^{-1}f_2f_3+f_4)=-\phi_p(-f_1^{-1}f_2f_3+f_4).
\end{eqnarray*}
As for $f(x) = (f_1 x+f_2 )^3$, $f_1\in\mathbb{F}_p^{\times}$, one has 
\begin{eqnarray*}
\sum_{x\in\Fp}\phi_p(f(x))&=&\sum_{x\in\Fp}\phi_p^3(f_1x+f_2)=\sum_{x\in\Fp}\phi_p(f_1x+f_2)=0.
\end{eqnarray*}
\end{Proof}
For $N \ge 2$, we set $\pi(N)$ to be the number of primes less than or equal to $N$. The following is Lemma 2.1 in \cite{Nagao}.
\begin{Lemma}
\label{Lem:convergence results}
Let $\{c_p\}_p$ be a bounded sequence of nonnegative numbers indexed by primes $p$. If one 
of the sequences $\left\{\displaystyle \frac{1}{\pi(N)}\sum_{\substack{p\le N\\ p:\textrm{ prime}}} c_p\right\}$ and $\left\{\displaystyle\frac{1}{N}\sum_{\substack{p\le N\\ p:\textrm{ prime}}} c_p\log p\right\}$ converges, then both converge to
the same limit.
\end{Lemma}

\begin{Lemma}
\label{Lem1}
Let $f(x)\in\Z[x]$ be a cubic polynomial. One has
$$\lim_{N\to\infty}\frac{1}{N}\sum_{\substack{p\le N\\ p:\textrm{ prime}}} \frac{\log p}{p}\sum_{x\in\Fp}\phi_p(f(x))=0.$$
\end{Lemma}
\begin{Proof}
We start with assuming that $f(x)$ has a repeated root over $\Q$. If $f(x)$ has a root of multiplicity three, then this means that the reduction of $f$ modulo every prime $p$ is of multiplicity three. Lemma \ref{Lem:symbols of cubics} implies that $\displaystyle \sum_{x\in\Fp}\phi_p(f(x))=0$, hence the result. 

If $f(x)$ has a root of multiplicity two, then this means that for all but finitely many primes the reduction of $f(x)$ has a double root modulo $p$. In view of 
Lemma \ref{Lem:symbols of cubics}, one has
$$\left|\sum_{\substack{p\le N\\ p:\textrm{ prime}}} \frac{\log p}{p}\sum_{x\in\Fp}\phi_p(f(x))\right|\le \sum_{\substack{p\le N\\ p:\textrm{ prime}}} \frac{\log p}{p}\left|\sum_{x\in\Fp}\phi_p(f(x))\right|\le\sum_{\substack{p\le N\\ p:\textrm{ prime}}} \frac{\log p}{p}<2+\log N,$$
where the last inequality is Mertens' first theorem. Thus the result follows in this case. 

Now we assume $f(x)$ does not have a multiple root over $\Q$. This means that for all but finitely
many primes, the reduction of $f(x)$ has no multiple roots. In particular, the equation $E:y^2 = f(x)$
represents an elliptic curve for all but finitely many primes. Setting $\Delta(E)$ to be the discriminant of $E$, one obtains
\begin{eqnarray*}
\left|\sum_{\substack{p\le N\\ p:\textrm{ prime}}} \frac{\log p}{p}\sum_{x\in\Fp}\phi_p(f(x))\right|&=&\left|- \sum_{\substack{p\le N\\ p\nmid \Delta(E) }} \frac{\log p}{p}a_p(E)+\sum_{\substack{p\le N\\ p\mid \Delta(E) }} \frac{\log p}{p}\sum_{x\in\Fp}\phi_p(f(x))\right| \\
&\le &\sum_{\substack{p\le N\\ p\nmid \Delta(E) }} \frac{\log p}{p}|a_p(E)|+\sum_{\substack{p\le N\\ p\mid \Delta(E) }} \frac{\log p}{p}\left|\sum_{x\in\Fp}\phi_p(f(x))\right| \\
 &\le & 2\sum_{\substack{p\le N\\ p\nmid\Delta(E)}} \frac{\log p}{\sqrt p}+\sum_{\substack{p\le N\\ p\mid \Delta(E) }} \frac{\log p}{p}\left|\sum_{x\in\Fp}\phi_p(f(x))\right|,
\end{eqnarray*}
where the last inequality follows from Hasse’s bound, and the second term in the third line is a constant. Now Lemma \ref{Lem:convergence results} shows that
$$\lim_{N\to\infty}\frac{1}{N}\sum_{\substack{p\le N\\ p:\textrm{ prime}}} \frac{\log p}{\sqrt p}=\lim_{N\to\infty}\frac{1}{\pi(N)}\sum_{\substack{p\le N\\ p:\textrm{ prime}}} \frac{1}{\sqrt p}.$$
One notices that
$$\sum_{\substack{p\le N\\ p:\textrm{ prime}}} \frac{1}{\sqrt p}\le \sum_{2\le i\le N} \frac{1}{\sqrt i}\le\int_1^N\frac{1}{\sqrt x}\, dx=2\sqrt{N}-1.$$
Finally, the Prime Number Theorem yields that
$$\lim_{N\to\infty}\frac{1}{N}\sum_{\substack{p\le N\\ p:\textrm{ prime}}} \frac{\log p}{\sqrt p}=0,$$
hence follows the result.
\end{Proof}

\section{The Mordell-Weil rank}
\label{sec:MordellWeil}
Let $\EE$ be a quadratic elliptic surface defined by $y^2 = a_3(T )x^3 + a_2(T )x^2 + a_1(T )x + a_0(T )$. One 
recalls that
$$a_3(T)x^3 + a_2(T)x^2 + a_1(T)x + a_0(T)=A(x)T^2+B(x)T+C(x)$$
where $A,B,C\in\Z[x]$, $\deg A=3$, and $\deg B,\deg C\le 3$. The following theorem allows one to evaluate $A_p(\EE)$. 

\begin{Theorem}
\label{Thm1}
Let $\EE$ be a quadratic elliptic surface defined over $\Q$ by the equation
$$y^2=a_3(T)x^3 + a_2(T)x^2 + a_1(T)x + a_0(T)=A(x)T^2+B(x)T+C(x)$$
where $A,B,C\in\Z[x]$, $\deg A=3$, and $\deg B,\deg C\le 3$. Then
\begin{eqnarray*}-pA_p(\EE)=p\sum_{\substack{x\in\Fp\\B(x)^2-4A(x)C(x)\equiv 0\textrm{ mod }p}}\phi_p(A(x))-\sum_{x\in\Fp}\phi_p(A(x))+R_p,
\end{eqnarray*}
where $|R_p|\le12.$
\end{Theorem}
\begin{Proof}
One starts with evaluating $A_p(\EE)$ recalling that if $p|\Delta(t)$, then $a_p(\EE_t)=0$:
\begin{eqnarray*}
-pA_p(\EE)&=&-\sum_{t\in\Fp}a_p(\EE_t)\\
&=&\sum_{\substack{t\in\Fp\\ p\nmid\Delta(t)}}\sum_{x\in\Fp}\phi_p(a_3(t)x^3+a_2(t)x^2+a_1(t)x+a_0(t))\\
&=&\sum_{x\in\Fp}\sum_{\substack{t\in\Fp\\ p\nmid\Delta(t)}}\phi_p(A(x)t^2+B(x)t+C(x))\\
&=&\left(\sum_{\substack{x\in\Fp\\ B(x)^2-4A(x)C(x)\equiv 0\textrm{ mod }p}}+\sum_{\substack{x\in\Fp\\ B(x)^2-4A(x)C(x)\not\equiv 0\textrm{ mod }p}}\right)\sum_{\substack{t\in\Fp\\ p\nmid\Delta(t)}}\phi_p(A(x)t^2+B(x)t+C(x)).
\end{eqnarray*}
Now one may use Lemma \ref{Lem:Legendre} to see that for a fixed $x\in\Fp$, one has
\begin{eqnarray*}
\sum_{t\in\Fp}\phi_p(A(x)t^2+B(x)t+C(x))=\begin{cases} 
      (p-1)\phi_p(A(x)) & \textrm{ if }B(x)^2-4A(x)C(x)\equiv 0\textrm{ mod }p \\
      -\phi_p(A(x)) &\textrm{ if }  B(x)^2-4A(x)C(x)\not\equiv 0\textrm{ mod }p.
   \end{cases}
\end{eqnarray*}
This yields that
\begin{eqnarray*}
-pA_p(\EE)&=&(p-1)\sum_{\substack{x\in\Fp\\ B(x)^2-4A(x)C(x)\equiv 0\textrm{ mod }p}}\phi_p(A(x))-\sum_{\substack{x\in\Fp\\ B(x)^2-4A(x)C(x)\not\equiv 0\textrm{ mod }p}}\phi_p(A(x))\\&-&\sum_{x\in\Fp}\sum_{\substack{t\in\Fp\\p|\Delta(t)}}\phi_p(A(x)t^2+B(x)t+C(x))
\end{eqnarray*}
Using Lemma \ref{Lem:symbols of cubics}, one obtains that
\begin{eqnarray*}
-pA_p(\EE)&=&p\sum_{\substack{x\in\Fp\\ B(x)^2-4A(x)C(x)\equiv 0\textrm{ mod }p}}\phi_p(A(x))-\sum_{x\in\Fp}\phi_p(A(x))-\sum_{x\in\Fp}\sum_{t\in\Fp}'\phi_p(A(x)t^2+B(x)t+C(x))
\end{eqnarray*}
where $\sum'$ means that the sum runs over $t\in \Fp$ such that $p$ is a prime of multiplicative reduction of $\EE_t$.  Again, according to Lemma \ref{Lem:symbols of cubics} and the fact that the discriminant of $\EE$ is a polynomial of degree $12$ in $T$, one knows that 
$$\left|\sum_{x\in\Fp}\sum_{t\in\Fp}'\phi_p(A(x)t^2+B(x)t+C(x))\right|=\left|\sum_{t\in\Fp}'\sum_{x\in\Fp}\phi_p(a_3(t)x^3+a_2(t)x^2+a_1(t)x+a_0(t))\right|\le 12,$$
hence the conclusion.
\end{Proof}

\begin{Remark}
\label{rem:1}
Let $\EE$ be a quadratic elliptic surface defined over $\Q$ by the following equation $$y^2=a_3(T)x^3 + a_2(T)x^2 + a_1(T)x + a_0(T)=A(x)T^2+B(x)T+C(x)$$
where $A,B,C\in\Z[x]$, $\deg A=3$, and $\deg B,\deg C\le 3$.  Then one has
$$|A_p(\EE)|\le \begin{cases}
6+\frac{2}{\sqrt p}+\frac{12}{p}& \textrm{if $A(x)$ has no repeated root}\\
6+\frac{13}{ p}& \textrm{if $A(x)$ has a root of multiplicity $2$}\\
6+\frac{12}{p}& \textrm{if $A(x)$ has a root of multiplicity $3$}
\end{cases}$$
This follows directly from Lemma \ref{Lem:symbols of cubics} and Theorem \ref{Thm1} together with Hasse’s bound when $f(x)$ has no repeated root.
\end{Remark}

Now the bounds of the Mordell-Weil rank of a quadratic elliptic surface follow.
\begin{Theorem}
\label{Thm:2}
Let $\EE$ be a quadratic elliptic surface defined over $\Q$ by the following equation $$y^2=a_3(T)x^3 + a_2(T)x^2 + a_1(T)x + a_0(T)=A(x)T^2+B(x)T+C(x)$$
where $A,B,C\in\Z[x]$, $\deg A=3$, and $\deg B,\deg C\le 3$.  We set
\[S_1=\{x\in\Q:B(x)^2-4A(x)C(x)=0\textrm{ and $A(x)$ is a nonzero rational square}\}.\]
The Mordell-Weil rank of $\EE(\Q(T))$ is given by 
\[\begin{cases} \rank(\EE(\Q(T)))=|S_1|& \textrm{ if } B(x)^2-4A(x)C(x)\textrm{ splits completely over } \Q\\
 |S_1|\le \rank(\EE(\Q(T)))\le |S_1|+\delta&\textrm{ otherwise}
\end{cases}\]
where $\delta$ is the number of irreducible factors of $B(x)^2-4A(x)C(x)$ over $\Q$ with degree at least $2$.

In particular,  $|S_1|\le \rank(\EE(\Q(T)))\le |S_1|+3.$
\end{Theorem}
\begin{Proof}
Since $\EE$ is a rational elliptic surface, it follows that $\displaystyle\rank(\EE(\Q(T)))=\lim_{X\to\infty}\frac{1}{X}\sum_{\substack{p\le X \\p:\textrm{ prime}}}-A_p(\EE)\log p$. As $\displaystyle\lim_{X\to\infty}\frac{1}{X}\sum_{\substack{p\le X \\p:\textrm{ prime}}}\sum_{x\in\Fp}\phi_p(A(x))\frac{\log p}{p}=0$ by Lemma \ref{Lem1}, Theorem \ref{Thm1} implies that
\[\rank(\EE(\Q(T)))=\lim_{X\to\infty}\frac{1}{X}\sum_{\substack{p\le X \\p:\textrm{ prime}}}\sum_{\substack{x\in\Fp\\ B(x)^2-4A(x)C(x)\equiv 0\textrm{ mod }p}}\phi_p(A(x))\log p.\]
Now we set \[S_2=\{x\in\Q:B(x)^2-4A(x)C(x)=0\textrm{ and $A(x)$ is not a nonzero rational square}\}.\]
Let $B(x)^2 - 4A(x)C(x) = h_1 (x)h_2 (x)g(x)$ where $h_i (x)\in \Z[x]$ is such that $h_i (x) = 0$ if and only if $x\in S_i$. It follows that $g(x)$ is a product of at most three irreducible polynomials each of degree
$\ge 2$. One has 
\[\sum_{\substack{x\in\Fp\\ B(x)^2-4A(x)C(x)\equiv 0\textrm{ mod }p}}\phi_p(A(x))=\sum_{x\in S_1}\phi_p(A(x))+\sum_{x\in S_2}\phi_p(A(x))+\sum_{\substack{x\in\Fp\\g(x)\equiv 0 \textrm{ mod }p}}\phi_p(A(x)).\]
One knows that for any $x \in S_1$, $\phi_p(A(x)) = 1$. Moreover, for any $x \in S_2$, either $A(x) = 0$, hence it does not contribute to the second sum, or $A(x)$ is not a square in $\Q$. In the latter case, one needs to evaluate the following limit for $x\in  S_2$
\begin{equation}\label{eq1}\lim_{X\to\infty}\frac{1}{X}\sum_{\substack{p\le X \\p:\textrm{ prime}}}\phi_p(A(x))\log p\le\lim_{X\to\infty}\frac{\log X}{X}\sum_{\substack{p\le X \\p:\textrm{ prime}}}\phi_p(A(x)). \end{equation}
In accordance to (5.79) of \cite{Iwaniec}, one has
\[\sum_{\substack{p\le X \\p:\textrm{ prime}}}\phi_p(A(x))=O(X(\log X)^{-\alpha})\]
for any $\alpha>0$ where the implied constant depends only on $\alpha$. In particular, by choosing $\alpha>1$, one obtains that the limit in (\ref{eq1}) is zero. Now, the latter argument together with the Prime Number Theorem and Lemma \ref{Lem:convergence results} yield that 
\[\rank(\EE(\Q(T)))=|S_1|+\lim_{X\to\infty}\frac{1}{X}\sum_{\substack{p\le X\\p:\textrm{ prime}}}\sum_{\substack{x\in\Fp\\g(x)\equiv 0 \textrm{ mod }p}}\phi_p(A(x))\log p.\]
Now, if $m(x)$ is an irreducible factor of $g(x)$,  $\deg m \ge 2$, the following inequality holds
$$\sum_{\substack{x\in\Fp\\m(x)\equiv 0 \textrm{ mod }p}}\phi_p(A(x))\le r_p(m(x))$$
where $r_p(m(x))$ is the number of zeros of $m(x)$ mod $p$.

Since $m(x)$ is irreducible, then by Chebotarev's density theorem, the density of primes $p$ such
that $m$ has $k$ roots mod $p$ is equal to the proportion of elements in the Galois group of $m(x)$ which fix $k$ roots when acting on the roots. As the Galois group acts transitively on the roots, the average number of fixed points is given by
\[\lim_{X\to\infty}\frac{1}{X}\sum_{\substack{p\le X \\p:\textrm{ prime}}}r_p(m(x))\log p=\lim_{X\to\infty}\frac{1}{\pi(X)}\sum_{\substack{p\le X \\p:\textrm{ prime}}}r_p(m(x))=1\]
where the first equality is due to Lemma \ref{Lem:convergence results}, whence the result follows.
\end{Proof}
\begin{Remark}
Let $\EE$ be a quadratic elliptic surface defined over $\Q$ by the equation $$y^2=a_3(T)x^3 + a_2(T)x^2 + a_1(T)x + a_0(T)=A(x)T^2+B(x)T+C(x)$$
where $A,B,C\in\Z[x]$, $\deg A=3$, and $\deg B,\deg C\le 3$.  
One has \[0\le \rank(\EE(\Q(T)))\le 6.\]
\end{Remark}
 If one wants to construct a quadratic elliptic surface with Mordell-Weil rank $r$, namely $1\le r\le 6$, then one considers the polynomial
$$D(x)=4(x-v_1)(x-v_2)(x-v_3)(x-v_4)(x-v_5)(x-v_6),\quad v_i\in\Z,\; v_i\ne v_j \textrm{ if }i\ne j.$$
One then chooses $A(x),B(x), C(x) \in \Z[x]$, with $\deg A(x)=3$ and $B(x), C(x)$ of degree at most $3$ such that $D(x) - B(x)^2 = -4A(x)C(x)$, $B(v_i)\ne 0$, $i=1,\cdots,r$; and $A(v_{i})$ is a nonzero rational sqaure for every $i=1,\ldots,r$. Now if we consider the elliptic surface
\[\EE:y^2=A(x)T^2+B(x)T+C(x)\]
over $\Q$, then according to Theorem \ref{Thm:2}, one has $\rank(\EE(\Q(T)))=r$. In \cite{Arms}, examples of such elliptic surfaces with maximal rank $6$ were given. In current work, we plan to present how the polynomials $A(x), B(x),$ and $C(x)$ may be chosen to produce an elliptic surface over $\Q$ with Mordell-Weil rank $r$, $1\le r\le 6$.


\begin{thebibliography}{MM}
\frenchspacing
\renewcommand{\baselinestretch}{1}


\bibitem{Arms} S. Arms, A. Lozano-Robledo and S. J. Miller, {\it Constructing one-parameter families of elliptic curves with moderate rank}, Journal of Number Theory, {\bf123} (2007), 388--402
\bibitem{Bat}
F.  Battistoni, S. Bettin and C. Delaunay,
 Ranks of elliptic curves over $\Q(T)$ of small degree in $T$, preprint

\bibitem{Fer1}
 S. Fermigier, {\it Un exemple de courbe elliptique d\'{e}finie sur $\Q$ de rang $\ge 19$}, 
 C. R. Acad. Sci. Paris S\'{e}r. 1, {\bf 315} (1992), 719--722

\bibitem{Fer}
S. Fermigier, 
{\it \'{E}tude exp\'{e}rimentale du rang de familles de courbes elliptiques sur $\Q$},
 Experimental Mathematics, {\bf 5} (1996), 119--130

\bibitem{Fer2}
S. Fermigier, 
{\it Une courbe elliptique d\'{e}finie sur $\Q$ de rang $\ge 22$}, Acta Arith., {\bf 82} (1997),  359--363

\bibitem{Hillgarter}
E. Hillgarter, {\it Rational points on conics}, Ph.D. thesis, Johannes Kepler Universit\"{a}t, Linz, Austria, 1996
\bibitem{Iwaniec} H. Iwaniec and E. Kowalski, {\it Analytic number theory}, Colloquium Publications, American Mathematical Society, Volume 53, 2004
\bibitem{Kollar} J. Koll\'{a}r and M. Mella, {\it Quadratic families of elliptic curves and unirationality of degree 1 conic bundles}, American Journal of Mathematics, {\bf 139} (2017), 915--936
\bibitem{Mackall}
B. Mackall, S. J. Miller, C. Rapti, and K. Winsor,
{\it Lower-order biases in elliptic curve Fourier coefficients in families}, in:
Frobenius Distributions: Lang-Trotter and Sato-Tate Conjectures (David Kohel and Igor Shparlinski, editors), Contemporary Mathematics, AMS, Providence, RI. {\bf 663} (2016)

\bibitem{Mest1}
J.-F. Mestre, {\it Courbes elliptiques de rang $\ge 11$ sur $\Q(t)$}, C. R. Acad. Sci. Paris S\'{e}r. 1, {\bf 313} (1991), 139--142
\bibitem{Mest2}
J.-F. Mestre, {\it Courbes elliptiques de rang $\ge 12$ sur $\Q(t)$}, C. R. Acad. Sci. Paris S\'{e}r. 1, {\bf 313} (1991), 171--174

\bibitem{Nagao1}
K. Nagao, 
{\it An example of elliptic curve over $Q$ with rank $\ge 20$}, Proc. Japan Acad.
Ser. A, {\bf 69} (1993), 291--293
\bibitem{Nagao2}
K. Nagao, {\it An example of elliptic curve over $\Q(T)$ with rank $\ge13$}, Proc. Japan Acad.
Ser. A, {\bf 70} (1994), 152--153
\bibitem{Nagao3}
K. Nagao and T. Kouya, {\it An example of elliptic curve over $Q$ with rank $\ge 21$},
Proc. Japan Acad. Ser. A, {\bf 70} (1994), 104--105

\bibitem{Nagao} K. Nagao, {\it $\Q(T)$-rank of elliptic curves and certain limit coming from the local points}, Manuscripta Mathematica, {\bf 92} (1997), 13--32

\bibitem{Rosen} M. Rosen and J. H. Silverman, {\it On the rank of an elliptic surface}, Invent.  Math. {\bf 133} (1998), 43--67


\end{thebibliography}
\end{document}